\documentclass[12pt,a4paper]{article}
\usepackage{wasysym}
\usepackage{a4wide}
\usepackage{epsfig}
\usepackage{subfigure}
\usepackage{color}
\usepackage{float}
\begin{document}

\begin{center}
{\large\bf On the random nature of

(prime) number distribution}

\medskip

$^{a,b}$Erika L. Alvarez\footnote{ealvarez@fyma.ucl.ac.be}, $^{b}$Jean Pestieau\footnote{pestieau@fyma.ucl.ac.be}

$^{a}$Instituto de F\' isica, Universidad Nacional Aut\' onoma de M\'exico,

Apartado postal 20364, 01000 M\'exico D.F., M\'exico

$^{b}$Institut de Physique Th\'eorique, Universit\'e catholique de Louvain,

Chemin du Cyclotron 2, B-1348 Louvain-la-Neuve, Belgique

\medskip

Preliminary version, 14/12/2004
\end{center}
\begin{center}Abstract\end{center}

Let $\pi(x)$ denote the number of primes smaller or equal to x. We compare $\sqrt{\pi}(x)$ with
$\sqrt{R}(x)$ and $\sqrt{\ell i}(x)$, where $R(x)$ and $\ell i(x)$ are the Riemann function and the logarithmic integral, respectively. We show a regularity in the distribution of the natural numbers in terms of a phase related to $(\sqrt{\pi}-\sqrt{R})$ and indicate how $\ell i(x)$ can cross $\pi(x)$ for the first time.

\section{Introduction}
\subsection{Preliminaries}
The function $\pi(x)$ is the function counting the number of primes smaller or equal to $x$. For example, $\pi(2)=1$, $\pi(3)=2$, $\pi(4)=2$, $\pi(5)=3$, \ldots In 1792, when he was 15 years old, Gauss proposed
$$\frac{x}{\ln x}$$
as an approximation to $\pi(x)$, which he refined afterwards \cite{gauss} to
$$\ell i(x)=PV\int_0^x \frac{dt}{\ln t}$$
where PV means the integral principal value. The function $\ell i (x)$ can also be written as $\ell i(x)=\int_{\mu}^x dt/\ln t$, with $\mu=1.4513692348\ldots$

 Later, Riemann \cite{riemann} improved the approximation with his Riemann function $R(x)$
defined as

$$R(x)=\sum_{n=1}^{\infty}\frac{\mu(n)}{n}\ell i(x^{1/n})$$
where $\mu$ is the M\"obius function \cite{mobius}, given by
$$\begin{array}{cccc}
 &  & 0 & \mbox{if } n\mbox{ has one or more primes repeated}\\
 \mu(n) & =  & 1 & \mbox{if }n=1 \\
 & & (-1)^k & \mbox{if }n\mbox{ is a product of } k \mbox{ different primes}
\end{array}$$
Riemann also proposed that \cite{zagier}
\begin{equation}\label{mera}
\pi(x)-R(x)=-\sum_{\rho}R(x^{\rho})
\end{equation}
where $\rho$ are the {\em trivial} and {\em non trivial} zeroes of the Riemann zeta function, $\zeta$, which is defined as
$$\zeta(s)=\sum_{k=1}^{\infty}\frac{1}{k^s}$$
for $\Re (s)>1$. Although Riemann did the analytical continuation of $\zeta$ to all the complex plane excepting the point $s=1$, an easier expression is given by \cite{hasse}
$$\zeta(s)=\frac{1}{1-2^{1-s}}\sum_{n=0}^{\infty}\frac{1}{2^{n+1}}\sum_{k=0}^n(-1)^k \frac{n!}{k!(n-k)!(k+1)^{s}}$$

The {\em trivial} zeroes of $\zeta$ are found easily from the relation \cite{hardy}
$$\zeta(1-s)=2(2\pi)^{-s}\cos\left(\frac{s\pi}{2}\right)\Gamma(s)\zeta(s)$$
because when $s=2n+1$, with $n$ an integer, $\zeta(-2n)=0$.

With respect to the {\em non trivial} zeroes, the Riemann hypothesis \cite{riemann} says that all
of them lie on the ``critical'' line, $\rho(t)=1/2+it$. It is one of the most important
problems of mathematics today.

The prime number theorem, proved independently by de la Vall\'ee-Poussin \cite{poussin} and Hadamard \cite{hadamard}, assures that

$$\lim_{x\to\infty}\,\frac{\pi(x)}{\ell i(x)}=\lim_{x\to\infty}\,\frac{\pi(x)}{R(x)}=\lim_{x\to\infty}\,\frac{\pi(x)\ln x}{x}= 1$$

Currently $\pi(x)$ has been computed up to $x\sim 10^{23}$. All the computed values of $\pi(x)$ today
satisfy the inequality $\ell i(x)>\pi(x)$. However, in 1914 Littlewood \cite{littlewood} showed that
 this inequality changes its sign infinitely often for very large $x$ \cite{carter}.

\subsection{Motivation}
\begin{figure}[H]
\centering{
\mbox{\subfigure{\epsfig{file=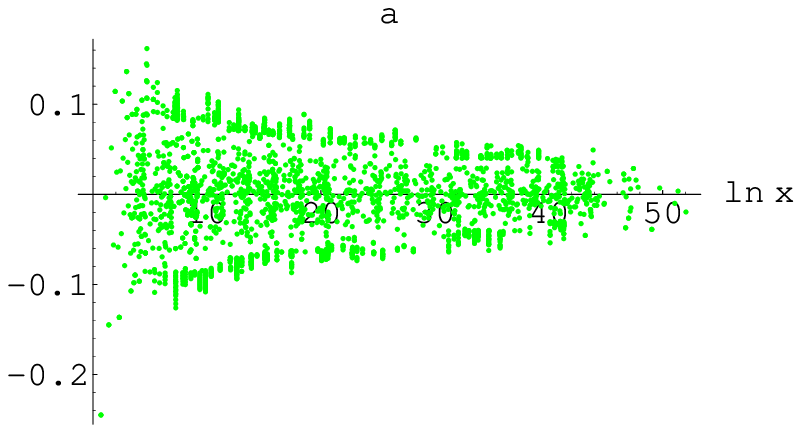,width=6.5cm,height=4.5cm}}\quad\subfigure{\epsfig{file=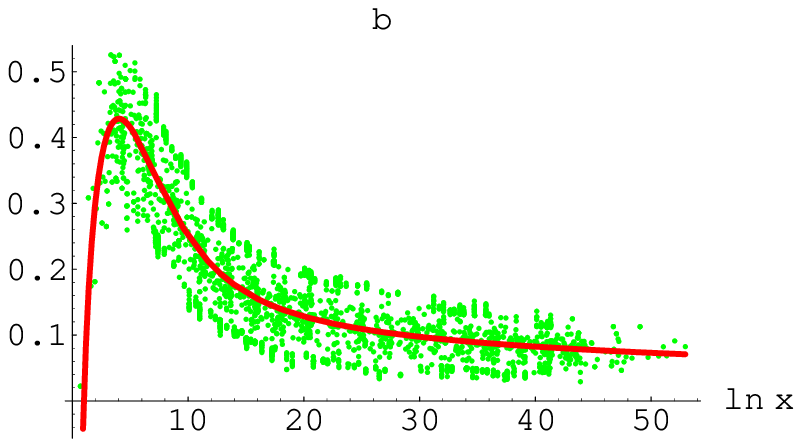,width=6.5cm,height=4.5cm}}}
}
\caption{(a) $\sqrt{\pi}(x)-\sqrt{R}(x)$ vs $\ln x$  and (b) $\sqrt{\ell i}(x)-\sqrt{\pi}(x)$ vs $\ln x$, in $x\in(2,10^{23})$, where the gross line is the function $\sqrt{\ell i}(x)-\sqrt{R}$(x) }
\label{piR}
\end{figure}

In general the absolute value of the difference between the function $\pi(x)$ and its approximations,
 $\ell i(x)$ or $R(x)$, although it is smaller than $\sim \sqrt{\pi}(x)$, is a number much greater
than the unity for large $x$. However, the absolute value of the difference between
 the square roots of $\pi(x)$ and of $\ell i(x)$ or between the square roots of $\pi(x)$ and of
$R(x)$ are smaller than $1$. Then these ones are what we will consider in order to have a better scope of
the approximations to $\pi(x)$. In Figure \ref{piR}.a, it is shown the difference
$\sqrt{\pi}(x)-\sqrt{R}(x)$ and the maximal difference between these functions is
$\sqrt{\pi}(2)-\sqrt R(2)=-0.244906$ when $x=2$. We see that $\sqrt{R}(x)$ averages very well
$\sqrt{\pi}(x)$. In Figure \ref{piR}.b it is shown the difference $\sqrt{\ell i}(x)-\sqrt{\pi}(x)$,
whose maximal height corresponds to the point $x=28$, where
$\sqrt{\ell i}(28)-\sqrt{\pi}(28)=0.525426$. The gross line represents the function
$\sqrt{\ell i}(x)-\sqrt{R}(x)$, which is the ``average'' of the points $\sqrt{\ell i}-\sqrt{\pi}$.
In both figures not all the points
are shown, there is a higher density in the center, a lot of external points are included to make
the border explicit. The points were calculated with Mathematica until $10^{12}$ and the rest were
taken from the tables of \cite{nicely}, which give values of $\pi(x)$ for numbers with three or four
significant digits, and so, the points shown in the border after $10^{12}$ are not necessarily the
points with the biggest difference $\mid\sqrt{\pi}-\sqrt{R}\mid$.

In section 2, our plan is to delimit the function $(\sqrt{\pi}-\sqrt{R})$ from above and below
with a tight function, in such a way that all the points remain inside the bounds, then,
to delimit the functions $(\sqrt{\ell i}-\sqrt{\pi})$ and $(\ell i -\pi)$, and finally to
discuss the statistical distribution of a phase defined in terms of the functions previously
mentioned.

\section{Discussion}
\subsection{$\sqrt{\pi}-\sqrt{R}$}
One can study the general characteristics of the function $\sqrt{\pi}(x)-\sqrt{R}(x)$.
The absolute value of this function is bounded with its maximal value
$\mid\sqrt{\pi}(2)-\sqrt R(2)\mid=0.244906$. So, we can propose that $\sqrt{\pi}(x)$ is given by
\begin{equation}\label{primera}
({\it i})\quad \sqrt{\pi}(x)=\sqrt{R}(x)+\eta(x)\cos\delta(x)\qquad \eta(x)>0
\end{equation}
where $\eta(x)$ is the envelope, and all the points of Figure \ref{piR}.a are delimited by this one.

\begin{figure}[H]
\centering{
\mbox{\subfigure{\epsfig{file=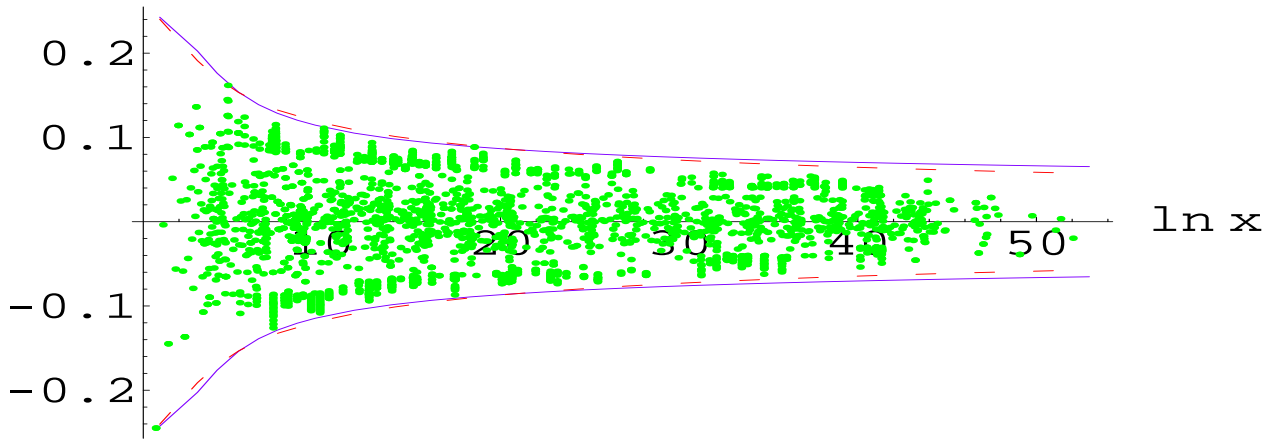,width=17cm,height=7.5cm}}}
}
\caption{ $\sqrt{\pi}(x)-\sqrt{R}(x)$ vs $\ln x$, envelopes $\eta(x)=0.2595/\ln\ln(x+15.9)$ (continuous line) and $\eta(x)=0.315647/[\ln(x+4.07206)]^{0.430202}$ (dashed line)}
\label{mano}
\end{figure}

\begin{figure}[H]
\centering{
\mbox{\subfigure{\subfigure{\epsfig{file=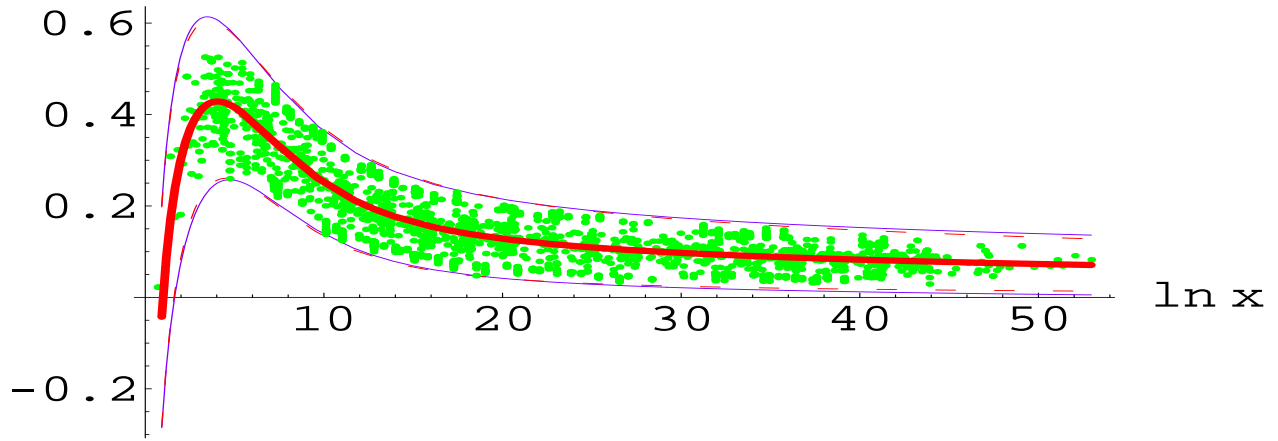,width=17cm,height=7.5cm}}}}
}
\caption{ $\sqrt{\ell} i(x)-\sqrt{\pi}(x)$ vs $\ln x$, envelopes $\eta(x)=0.2595/\ln\ln(x+15.9)$ (continuous line) and $\eta(x)=0.315647/[\ln(x+4.07206)]^{0.430202}$ (dashed line), and the function $\sqrt{\ell i}(x)-\sqrt{R}(x)$ (gross line)}
\label{mano2}
\end{figure}

Other parameterization is
\begin{equation}\label{segunda}
({\it ii})\quad a(x)=\sqrt{R}(x)+\eta(x)e^{i\delta(x)}\qquad \eta(x)>0,\quad\mid a(x)\mid^2=\pi(x)
\end{equation}
this last one puts in evidence the parameterization in terms of an amplitude $\eta(x)$ and a
phase $\delta(x)$. Equation (\ref{segunda}) implies
\begin{equation}\label{tercera}
\pi(x)=R(x)+2\eta(x)\cos\delta\sqrt{R}(x)+\eta^2(x)
\end{equation}
Observe that, when $\delta(x)=0$ or $\pi$, Equations (\ref{primera}) and (\ref{segunda}) coincide.
The first proposal for $\eta(x)$ is the function
\begin{equation}\label{loglog}
\eta_1(x)=\frac{0.2595}{\ln\ln(x+15.9)}
\end{equation}
However, from the work of \cite{carter} we know that the first zero of the function
$\sqrt{\ell i}(x)-\sqrt{\pi}(x)$ happens before $x=1.3982\times 10^{316}$, and may be much earlier.
The function of Equation (\ref{loglog}) crosses x axis around  $x=10^{65}$. A function that crosses
x axis around $x=1.3982\times 10^{316}$, is
\begin{equation}\label{logpotencia}
\eta_2(x)=\frac{0.315647}{[\ln(x+4.07206)]^{0.430202}}
\end{equation}
If  $\sqrt{\ell i}(x)-\sqrt{\pi}(x)$ crossed the axis before, $\eta(x)$ would be a function between
the ones defined in Equation (\ref{loglog}) and Equation (\ref{logpotencia}).
In Figure \ref{mano} it is shown the points $\sqrt{\pi}(x)-\sqrt{R}(x)$ with the two bounds and in Figure \ref{mano2} it is shown the points $\sqrt{\ell i}(x)-\sqrt{\pi}(x)$, with its ``average'' function  $\sqrt{\ell i}(x)-\sqrt{R}(x)$, where the borders are given by
$$(\sqrt{\ell i}-\sqrt{\pi})_{\mbox{{\scriptsize max,min}}}=\sqrt{\ell i}-\sqrt{R}\pm\eta$$

\subsection{$\ell i-\pi$}
We can delimit $\ell i-\pi$ from above and below.

From Equation (\ref{primera}) and Equation (\ref{tercera}) and using $\ell i-\pi=\ell i-R+R-\pi$ one has that
\begin{equation}\label{exacto}
\ell i-(\sqrt{R}+\eta)^2\leq \ell i -\pi\leq \ell i-(\sqrt{R}-\eta)^2
\end{equation}
Using the fact that in the limit of large $x$, $\ell i(x)-R(x)\to \sqrt{x}/(\ln x)$, $\sqrt{R}\approx \sqrt{x/\ln x}$ and that $\eta^2$ is negligible, one has

\begin{equation}\label{aproximado}
\frac{\sqrt{x}}{\ln x}-2\eta\sqrt{\frac{x}{\ln x}}<\ell i-\pi<\frac{\sqrt{x}}{\ln x}+2\eta\sqrt{\frac{x}{\ln x}}
\end{equation}
and then, if there are values where $\ell i (x)$ is smaller than $\pi(x)$, then $\eta(x)$ must decrease in a slower way than $1/(2\sqrt{\ln x})$, as it happens with Equation (\ref{loglog}) and Equation (\ref{logpotencia}).

\begin{figure}[H]
\centering{
\mbox{\subfigure{\epsfig{file=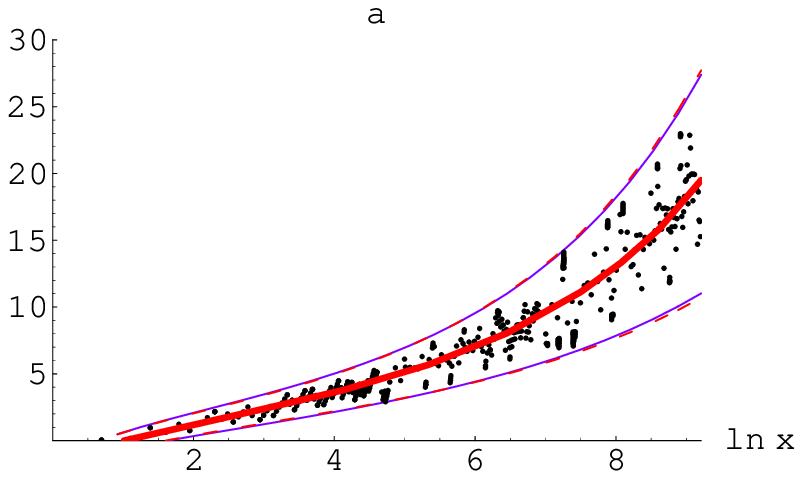,width=5cm,height=4.5cm}}
\subfigure{\epsfig{file=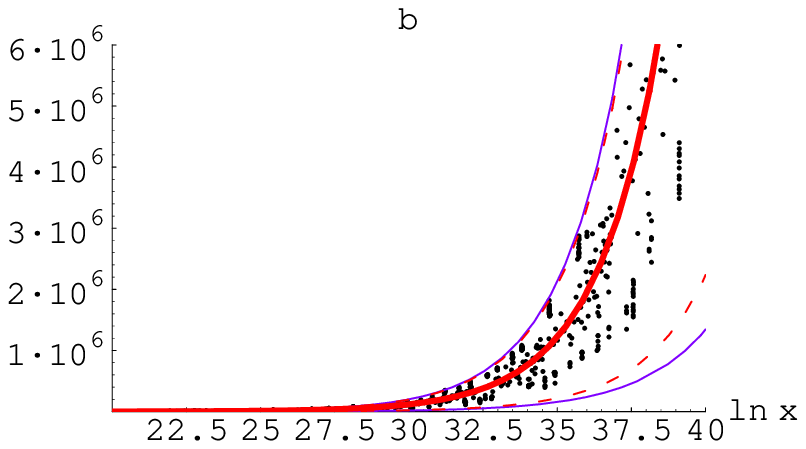,width=5cm,height=4.5cm}}
\subfigure{\epsfig{file=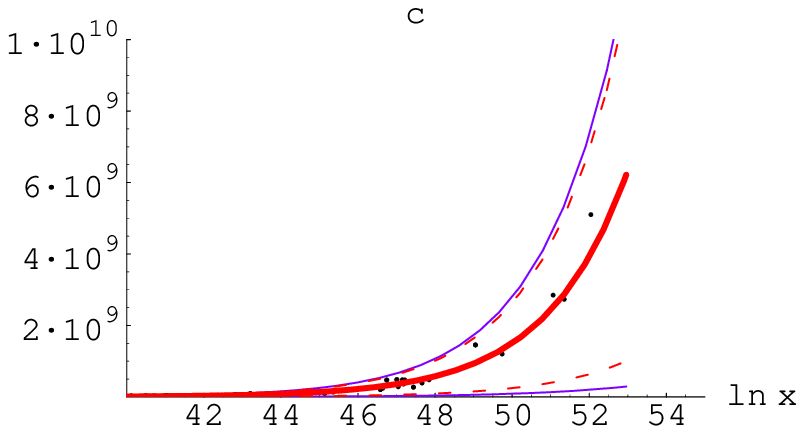,width=5cm,height=4.5cm}}}
}
\caption{ $\ell i(x)-\pi(x)$ vs $\ln x$, in (a) $x\in(2,10^{4})$, (b) $x\in(5\times 10^8,2\times 10^{17})$ and (c) $x\in(2\times 10^{17}, 8\times 10^{23})$, with $\eta(x)=0.2595/\ln\ln(x+15.9)$ (continuous line), $\eta(x)=0.315647/[\ln(x+4.07206)]^{0.430202}$ (dashed line) y $\ell i(x)-R(x)$ (gross line)  }
\label{gatos}
\end{figure}

 In Figures \ref{gatos}.b and \ref{gatos}.c it is shown $\ell i(x)-\pi(x)$ using for their
bounds Equation (\ref{aproximado}), with $\eta(x)$ given by Equation (\ref{loglog}) and Equation (\ref{logpotencia}). The bounds of Equation (\ref{aproximado}) only work for large $x$, when $R(x)\approx \ell i(x)-(1/2)\ell i(x^{1/2})$. For small $x$, Equation (\ref{aproximado}) is not valid, and we use directly the bounds (\ref{exacto}), and in Figure \ref{gatos}.a we show the later ones in the interval $x\in(2,10^4)$. The gross line corresponds to the ``average'' function $(\ell i(x)-R(x))$.

\subsection{$\cos\delta$}
With a sample of the first natural numbers one averages the functions $\sqrt{\pi}-\sqrt{R}$ and $\pi-R$.
The values of Table \ref{tbuno} are obtained for different sample sizes. In this table, $\sigma(f)$
is the standard deviation, $\sigma\equiv\sqrt{\langle f^2\rangle-\langle f\rangle^2}$, with $f$
equal to $(\sqrt{\pi}-\sqrt{R})$ or to $(\pi-R)$. We see that
$\langle{\tiny\sqrt{\pi}\!-\!\sqrt{R}}\rangle$ is a small number bigger than zero and has a small
variation in the different intervals.

 Working out the value of $\cos\delta$ in both cases, Equations (\ref{primera}) and (\ref{tercera}), one has
$$\cos\delta=\frac{\sqrt{\pi}(x)-\sqrt{R}(x)}{\eta(x)}\quad \mbox{y}\quad\cos\overline{\delta}=\frac{\pi(x)-R(x)-\eta^2(x)}{2\sqrt{R}(x)\eta(x)}$$
respectively and taking the first definition of $\eta(x)\equiv \eta_1(x)$, Equation (\ref{loglog}), one has
the averages of Table \ref{tbdos} in the intervals $x\in(2,100)$,\ldots, $x\in(2,10^6)$.

The results of Table \ref{tbdos} show that the averages remain approximately constant. With respect
 to the width of $\sigma$ of the distribution, as to the average of $\cos\delta$ absolute value,
 the difference in the
parameterizations of Equation (\ref{primera}) and Equation (\ref{segunda}) is negligible. Also, although
for the first intervals the difference in the average $\langle\cos\delta\rangle$ is bigger, as $x$ grows
the averages in the two parameterizations get closer,
because in general the ratio $\eta^2/|\pi-R|\ll 1$. From now on, we will keep the parameterization of
Equation (\ref{primera}).

Taking the other proposal of $\eta(x)\equiv\eta_2(x)$, Equation (\ref{logpotencia}), the averages of
 Table \ref{tbtres} are found. In this table, the average value of $\cos\delta$ is not very
different from the previous parameterization, being consistent with a small positive number.

In order to see the weight of the different sets of numbers with respect to $\cos\delta$, in Table
\ref{todos} we give the average of $\cos\delta$ for natural, prime, even and odd (without
primes) numbers. We see that as $x$ grows, the prime distribution, which has a higher $\cos\delta$
average, has a smaller weight, because
the ratio of prime to natural numbers decreases approximately as $\pi(x)/x\sim1/\ln x$. So,
 the average of $\cos\delta$ for the even and odd natural numbers will be approximately the same
for large $x$.

Let us take $\eta(x)$ given by Equation (\ref{loglog}): if we divide $\cos\delta$ in the intervals
\linebreak $(-1,-0.95),(-0.95,-0.85),\ldots,(0.85,0.95),(0.95,1)$, we find distributions of
Table \ref{reparto}. They give the number of positive integers whose $\cos\delta$ falls in one of
these intervals, we count them in $4$ different sample sizes: $(2,10^3)$, $(2,10^4)$, $(2,10^5)$ and
 $(2,10^3)$.

In Figures \ref{mano3} and \ref{mano4} it is shown distributions of $\cos\delta$ as
explained in the previous paragraph. We have normalized them to have the total area of the bars
equal to one. For example, for the natural numbers between $(2,10^3)$, there are $47$ numbers whose
$\cos\delta$ falls in the interval $(-0.45,-0.35)$. We divide these $47$ numbers by the sample
total number, $999$, to obtain the relative frequency and  multiply by $10$, because
the size of each interval is $0.1$ (except for the intervals $(-1,-0.95)$ and $(0.95,1)$).

The distribution is gaussian, and from Table 2,
 the width appears to have the same value, $\sigma=0.28$, it does not matter the number of positive integers
with which we take the average. The average seems to
stabilize around $\langle\cos\delta\rangle=0.014$. In all the figures, Figures \ref{mano3} and
\ref{mano4}, we used the same Gaussian
with width $\sigma=0.28$, average $\langle\cos\delta\rangle=0.014$ and height
$1/(\sqrt{2 \pi}\sigma)=1.425$, and the fit of the Gaussian is in a very good agreement with the data.

Finally, from Equation (\ref{primera}) and Equation (\ref{tercera})
$$\pi-R\approx 2\sqrt{R}\,\eta\cos\delta$$
then $(\pi(x)-R(x))/( 2\sqrt{R}(x)\eta(x))$ follows the same Gaussian distribution.

\section{Conclusions}
\begin{figure}[H]
\centering{
\mbox{\subfigure{\subfigure{\epsfig{file=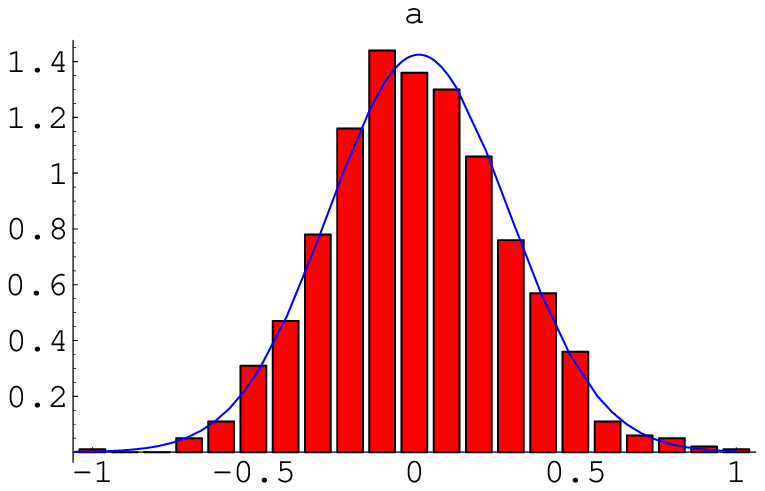,width=7cm,height=5.5cm}}\quad\subfigure{\epsfig{file=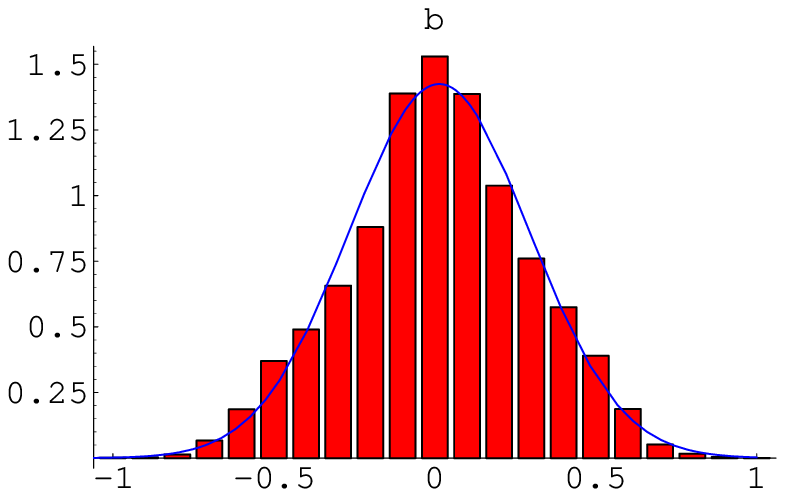,width=7cm,height=5.5cm}}}}
}
\caption{distribution of $\cos\delta$ with $\eta(x)=0.2595/\ln\ln(x+15.9)$, where the relative frequence of $\cos\delta$ has been counted in the intervals $(-1,0.95),(-0.95,-0.85)\ldots$, the Gaussian is represented by the continuous line,
(a) for the first $10^3$ natural numbers (except $1$) and (b) for the first $10^4$ ones  }
\label{mano3}
\end{figure}

\begin{figure}[H]
\centering{
\mbox{\subfigure{\subfigure{\epsfig{file=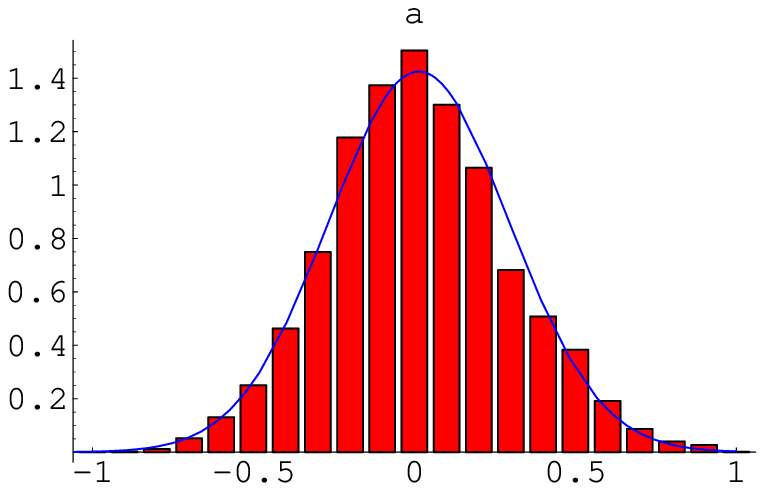,width=7cm,height=5.5cm}}\quad\subfigure{\epsfig{file=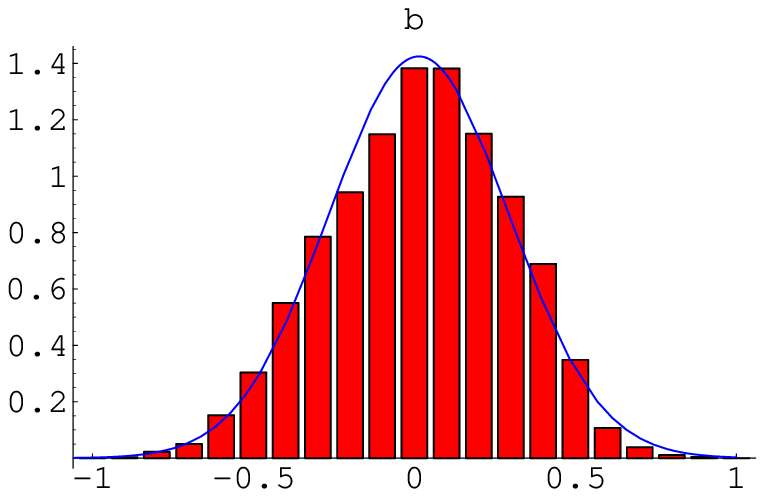,width=7cm,height=5.5cm}}}}
}
\caption{ the same as in the figure (\ref{mano3}), (a) for the first $10^5$ natural numbers
(except $1$) and (b) for the first $10^6$ ones }
\label{mano4}
\end{figure}

With two parameters, one amplitud $\eta$ and a phase $\delta$, we study the properties of the roots of the
functions $\pi$, $\ell i$ and $R$, using Equations (\ref{primera}) and (\ref{segunda}). With $\eta$, we delimit the differences ($\sqrt{\pi}-\sqrt{R}$), $(\sqrt{\ell i}-\sqrt{\pi})$ and $(\ell i-\pi)$. Concerning the last one, we know from the data that
$(\ell i-\pi<\sqrt{\pi})$, and in Equation (\ref{exacto}) we give a more precise relation. We find that
 $\cos\delta$, follows a Gaussian distribution, that shows a stable random behavior of the function $\pi(x)$,
 see Figures \ref{mano3} and \ref{mano4}. Taking different sample sizes, $\cos\delta$ distribution remains
constant. The question is if the Gaussian shape remains constant as  $x$ grows.

\vspace{1cm}

{\Large\bf\noindent Appendix }

\vspace{0.7cm}

To see how the natural numbers accommodate in the different $\cos\delta$ intervals, we give as an example
the first hundred in Table \ref{cien}, where the prime numbers have been underlined.

\begin{figure}[H]
\centering{
\mbox{\subfigure{\epsfig{file=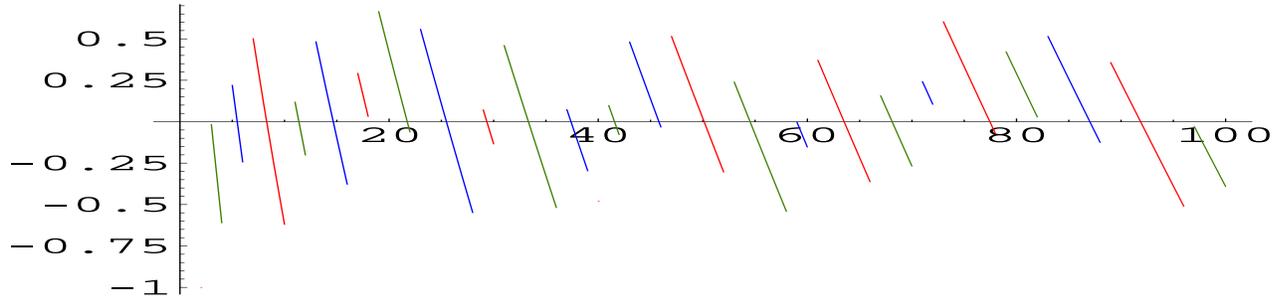,width=17cm,height=4.5cm}}}
}
\caption{  $\cos\delta=(\sqrt{\pi}(x)-\sqrt{R}(x))/\eta(x)$ vs $x$, $x\in(2,100)$, $\eta(x)=0.2595/\ln\ln(x+15.9)$ }
\label{luego}
\end{figure}

We can see that, each time there is a new prime number, $\cos\delta$ increases, and meanwhile $\pi(x)$ remains constant, until the next prime number, the following integers accommodate in intervals with smaller $\cos\delta$. So, Table 5 and
 Figures \ref{mano3} and \ref{mano4} show that the way of appearance of the prime numbers
implies the randomness of the natural numbers with respect to $\cos\delta$. In Figure \ref{luego},
we give a pictorial representation of how the first one hundred natural numbers (except 1) are accomodated, where the lines join points with the same $\pi(x)$

\begin{figure}[H]
\centering{
\mbox{\subfigure{\epsfig{file=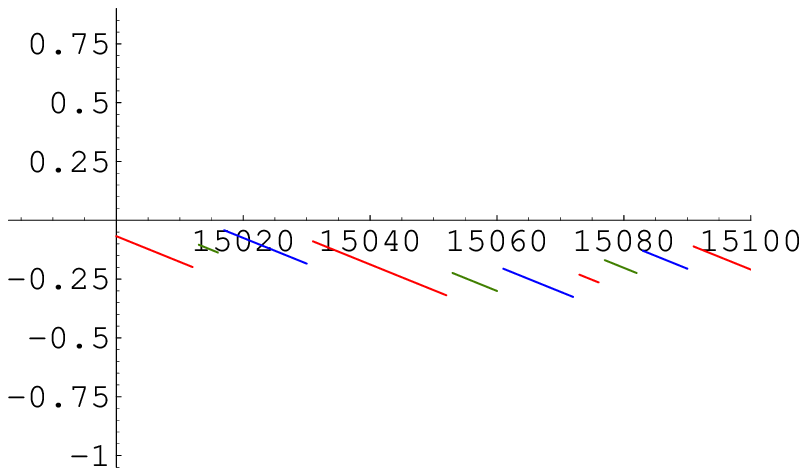,width=17cm,height=4.5cm}}}
}
\caption{  $\cos\delta=(\sqrt{\pi}(x)-\sqrt{R}(x))/\eta(x)$ vs $x$, $x\in(15\,000,15\,100)$, $\eta(x)=0.2595/\ln\ln(x+15.9)$ }
\label{manana}
\end{figure}

\begin{figure}[H]
\centering{
\mbox{\subfigure{\epsfig{file=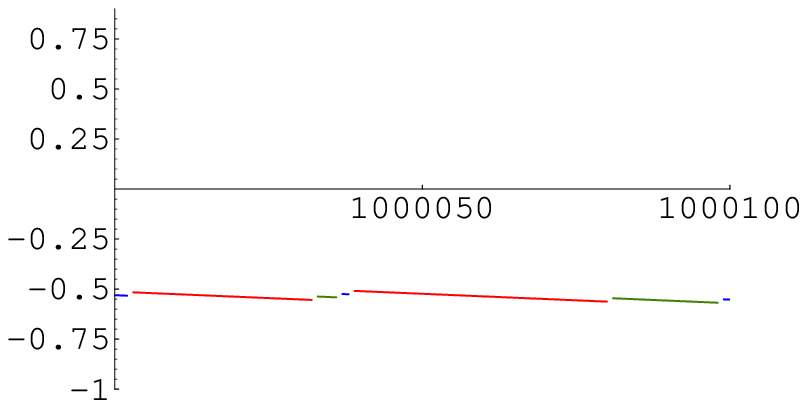,width=17cm,height=4.5cm}}}
}
\caption{  $\cos\delta=(\sqrt{\pi}(x)-\sqrt{R}(x))/\eta(x)$ vs $x$, $x\in(1\,000\,000,1\,000\,100)$, $\eta(x)=0.2595/\ln\ln(x+15.9)$ }
\label{manana2}
\end{figure}

That $\cos\delta$ decreases each time  $\pi(x)$ remains constant, while a new prime number does not appear,
it is because the function $\sqrt{R}(x)$ is a monotone growing function. With the appearance of the new prime number, $\cos\delta$ increases and the cycle is repeated. The rate with which $\cos\delta$ decreases is given by its derivative, and as the derivative of $R(x)$ is
$$ \frac{dR}{dx}=\sum_{n=1}^{\infty}\frac{\mu(n)}{n x^{(n-1)/n}\ln x}= \frac{1}{\ln x}\left(1-\frac{1}{2 x^{1/2}}-\frac{1}{3 x^{2/3}}-\ldots\right)$$
then, with the parameterization of Equation (\ref{loglog}) the derivative of $\cos\delta(x)$ for a constant $\pi(x)$ is
$$\frac{d\cos\delta}{dx}=\frac{1}{0.2561}\left\{\frac{(\sqrt\pi-\sqrt R)}{(x+15.5)\ln(x+15.5)}-\frac{\ln\ln(x+15.5)}{2\sqrt R\ln x}\left(1-\frac{1}{2x^{1/2}}-\ldots\right)\right\}$$
while for the parameterization of Equation (\ref{logpotencia}) is
$$\frac{d\cos\delta}{dx}=\frac{[\ln(x+4.07)]^{0.43}}{0.3156}\left\{0.43\frac{(\sqrt\pi-\sqrt R)}{(x+4.07)\ln(x+4.07)}-\frac{1}{2\sqrt R\ln x}\left(1-\frac{1}{2x^{1/2}}-\ldots\right)\right\}$$
in both cases the derivative is dominated by the negative term, as it is expected, and decreases in
absolute value as $x$ grows. In Figures \ref{manana} and \ref{manana2}, some other intervals of $100$ numbers are compared for larger $x$ where it is seen that $\cos\delta$ gets more horizontal, this is because there is a bigger number of points with the same $\pi(x)$, also, although at the beginning there are ``jumps'' when one goes from $\pi(p-1)$ to $\pi(p)$, whose
difference is one, as $x$ increases, $\cos\delta$ turns into a softer function, because $(\pi(p)-\pi(p-1))/\pi(p)\to 0$.

\medskip

We would like to thank Laurent Jacques, Gabriel L\'opez Castro,
Thomas R. Nicely, Mat\' ias Moreno and V\'{\i}ctor Romero for his kind
assistance in this work. We would also like to acknowledge the
support of project IN-120602 of the Direcci\'{o}n de Asuntos del
Personal Acad\'{e}mico of the Universidad Nacional Aut\'{o}noma de M\'{e}xico.

\pagebreak

\begin{center}{\large \bf Tables}\end{center}

\begin{center}
\begin{table}[H]
$$\begin{array}{|c||c|c||c|c|}\hline & \langle{\tiny\sqrt{\pi}\!-\!\sqrt{R}}\rangle &  \sigma(\sqrt{\pi}\!-\!\sqrt{R}) &\langle\pi-R\rangle &   \sigma(\pi\!-\!R)  \\\hline 2-10^2 & 0.001889   & 0.062256 &  0.033137  & 0.403334\\2-10^3 & 0.001363 &   0.042803 &  0.013466  & 0.714523 \\2- 10^4 & 0.001302  & 0.035624 & 0.050812  & 1.72635 \\2- 10^5 & 0.001529  & 0.031321 & 0.25608  & 4.23254 \\2- 10^6 & 0.001405 & 0.028509 &  0.705741  & 11.1907\\ \hline
\end{array}$$
\caption{averages $\langle{\tiny\sqrt{\pi}\!-\!\sqrt{R}}\rangle$ and $\langle\pi-R\rangle$ in $5$ intervals, $\sigma$ is the standard deviation}
\label{tbuno}
\end{table}
\end{center}

\begin{center}
\begin{table}[H]
$$\begin{array}{|c||c|c|c||c|c|c|c|}\hline  & \langle\cos\delta\rangle\!=\!\langle\frac{\sqrt{\pi}-\sqrt{R}}{\eta_1}\rangle & \sigma  &\langle\mid\cos\delta\mid\rangle\! & \langle\cos\overline{\delta}\rangle\!=\!\langle\frac{\pi\!-\!R-\!\eta_1^2}{2\sqrt{R}\eta_1}\rangle & \overline{\sigma}  &\langle\mid\cos\overline{\delta}\mid\rangle\! \\\hline 2-10^2   &  0.014402 & 0.315325   & 0.254145 & -0.010719 & 0.315370 & 0.253362 \\2-10^3 & 0.008304 & 0.280109 & 0.223332 & -0.000662 & 0.280161 & 0.223363\\2- 10^4 &  0.009965 & 0.283603 & 0.224534 & 0.006999 & 0.282839  & 0.224425 \\ 2-10^5  & 0.014043  & 0.281287 & 0.222306 & 0.013073 & 0.281312 & 0.222302 \\ 2-10^6   & 0.014057 & 0.278975  &0.227005  & 0.013740 & 0.278979 & 0.226989  \\\hline\end{array}$$
\caption{averages of $\cos\delta$ defined by Equation (\ref{primera}) and of $\cos\overline{\delta}$
 given by Equation (\ref{tercera}), where $\eta_1$ is given by Equation (\ref{loglog})}
\label{tbdos}
\end{table}
\end{center}

\begin{center}
\begin{table}[H]
 $$\begin{array}{|c||c|c|c|}\hline   & \langle\cos\delta\rangle\!=\!\langle\frac{\sqrt{\pi}-\sqrt{R}}{\eta_2}\rangle & \sigma  &\langle\mid\cos\delta\mid\rangle\! \\\hline 2-10^2 & 0.015587 & 0.327365 & 0.264286  \\ 2-10^3 & 0.008606 & 0.279154 & 0.222093 \\ 2- 10^4 & 0.009728 & 0.274722  & 0.217424 \\2- 10^5 & 0.013529 & 0.270749 &0.213973 \\2- 10^6 & 0.013608 & 0.269512 & 0.219318 \\ \hline\end{array}$$
\caption{averages of $\cos\delta$ defined in Equation (\ref{primera}), with $\eta_2$ given in
Equation (\ref{logpotencia})}
\label{tbtres}
\end{table}
\end{center}

\begin{center}
\begin{table}[H]
$$\begin{array}{|c|c|cc|c|cc|c|}
\hline
 & \mbox{all} &\multicolumn{2}{c|}{\mbox{primes}}  &\mbox{even} &
\multicolumn{2}{c|}{\mbox{odd without }1}    &\mbox{odd} \\
 &\mbox{without } 1 & \multicolumn{2}{c|}{}          & \mbox{without }2 &
\multicolumn{2}{c|}{\mbox{and without primes}} &\mbox{without }1 \\ \hline
2-10^2 & 0.014402 &0.256581 & (25)&-0.073010  & -0.056451 & (25)& 0.122513  \\
2-10^3 & 0.008304 &0.182444 & (168)&-0.027533 & -0.025953 & (332)&0.046160 \\
2-10^4 & 0.009965 &0.099122&(1229)&-0.002571 & -0.002474 & (3771)& 0.022703 \\
2-10^5 &0.014043  &0.051776 & (9592)& 0.009936 & 0.010169 &(40408)&0.018171\\
2-10^6 &0.014057 &0.028670 &(78498)& 0.012749 & 0.012886 &(421502)& 0.015366 \\
\hline
\end{array}$$
\caption{average of $\cos\delta$, according to the set of positive integers under which the
average is taken, the numbers in parentheses are the number of positive integers in the given set}
\label{todos}
\end{table}
\end{center}

\begin{center}
\begin{table}[H]
$$\begin{array}{|c||c|c|c|c|}\hline \cos\delta & 2\!-\!10^3 & 2\!-\!10^4 &  2\!-\!10^5 &  2\!-\!10^6 \\\hline  (-1,-0.95) & 1 & 2 & 3 & 3\\(-0.95,-0.85) & 0 & 3 & 26 & 331 \\(-0.85,-0.75) & 0 & 14 & 120 & 2230\\ (-0.75,-0.65) & 5 & 67 & 522 & 5024 \\ (-0.65,-0.55) & 11 & 186 & 1303 & 15247 \\ (-0.55,-0.45) & 31 & 370 & 2504 & 30391\\ (-0.45,-0.35) & 47 & 490 & 4630 & 55051\\ (-0.35,-0.25) & 78 & 657 & 7490 & 78559\\ (-0.25,-0.15) & 116 & 880 & 11776 & 94341\\ (-0.15,-0.05) & 144 & 1389 & 13740 & 114888 \\ (-0.05,0.05) & 136 & 1530 & 15040 & 138262\\ (0.05,0.15) & 130 & 1387 & 13006 & 138171\\ (0.15,0.25) & 106 & 1038 & 10645 & 115027\\ (0.25,0.35) & 76 & 760 & 6816 & 92749\\ (0.35,0.45) & 57 & 575 & 5082 & 68886\\(0.45,0.55) & 36 & 390 & 3835 & 34856 \\ (0.55,0.65) & 11 & 187 & 1915 & 10700\\ (0.65,0.75) & 6 & 52 & 871 & 3833\\(0.75,0.85) & 5 & 17 & 397 & 1086 \\(0.85,0.95) & 2 & 4 & 267 & 373\\(0.95, 1) & 1 & 1 & 11 &11\\\hline\end{array}$$
\caption{number of positive integers with $\cos\delta$ in the intervals $(-1,-.95),(-0.95,-0.85),\ldots$ for differents samples: from $(2,10^3)$ to $(2,10^6)$}
\label{reparto}
\end{table}
\end{center}

\begin{center}
\begin{table}[H]
$$\begin{array}{c|l}\cos\delta & \\ (-1,-0.95) & {\bf\underline{2}} \\ (-0.65,-0.55) & 4,10 \\
(-0.55,-0.45) & 28,36,40,58,96 \\ (-0.45,-0.35) & 16,57,66,95,100 \\ (-0.35,-0.25) & 9,27,35,39,52,70,94,99
\\ (-0.25,-0.15) & 6,12,56,60,65,98 \\
(-0.15,-0.05) & 15,22,26,30,34,38,42,51,55,64,69,78,88,93 \\
(-0.05,0.05) & {\bf\underline{3}},18,46,50,{\bf\underline{59}},68,82,87,92,{\bf\underline{97}} \\
(0.05,0.15) & 8,{\bf\underline{11}},25,{\bf\underline{29}},33,{\bf\underline{37}},{\bf\underline{41}},
45,54,63,72,77,86,91 \\
(0.15,0.25) &{\bf\underline{5}},14,21,49,{\bf\underline{53}},62,{\bf\underline{67}},{\bf\underline{71}},
76,81,90 \\
(0.25,0.35) & {\bf\underline{17}},24,32,44,48,75,80,85 \\
(0.35,0.45) & 20,{\bf\underline{61}},{\bf\underline{79}},84,{\bf\underline{89}} \\
(0.45,0.55) & {\bf\underline{7}},{\bf\underline{13}},{\bf\underline{31}},{\bf\underline{43}},
{\bf\underline{47}},74,{\bf\underline{83}} \\
(0.55,0.65) & {\bf\underline{23}},\bf{\underline{73}} \\ (0.65,0.75) & {\bf\underline{19}}\\\end{array}$$
\caption{distribution of the first one hundred natural numbers (without $1$) in the different
intervals of $\cos\delta$, the primes are underlined.}
\label{cien}
\end{table}
\end{center}

\pagebreak

\begin{center}{\large\bf Figures}\end{center}

\end{document}